\documentclass[sn-mathphys-num]{sn-jnl}% Math and Physical Sciences Numbered Reference 

\usepackage{indentfirst,enumerate,amssymb,amsfonts,amsmath,amsthm,mathrsfs,dsfont}
\usepackage{color}
\usepackage{graphicx}%
\usepackage{xcolor}%

\numberwithin{equation}{section}

\theoremstyle{thmstyleone}
\newtheorem{theorem}{Theorem}[section]% meant for sectionwise numbers
\newtheorem{proposition}[theorem]{Proposition}
\newtheorem{corollary}[theorem]{Corollary}

\theoremstyle{thmstyletwo}
\newtheorem{example}{Example}[section]
\newtheorem{remark}{Remark}

\theoremstyle{thmstylethree}
\newtheorem{definition}{Definition}[section]
\newcommand{\C}{\mathbb{C}}
\newcommand{\TT}{{\mathbb{T}}}

\newcommand{\TrSs}{\TT^r\times\mathbb{S}^{3s}}

\def\jp#1{{\left\langle{#1}\right\rangle}}

\begin{document}

	\title[Systems on compact Lie groups]{Diagonal systems of differential operators  on compact Lie groups}


\author[1]{\fnm{P. L.} \sur{Dattori da Silva}}\email{dattori@icmc.usp.br}
\equalcont{These authors contributed equally to this work.}

\author[2]{\fnm{A.} \sur{Kirilov}}\email{akirilov@ufpr.br}
\equalcont{These authors contributed equally to this work.}

\author*[1]{\fnm{R.} \sur{Paleari da Silva}}\email{ricardopaleari@gmail.com}
\equalcont{These authors contributed equally to this work.}

\affil*[1]{\orgdiv{Departamento de Matem\'atica}, \orgname{Instituto de Ci\^encias Matem\'aticas e de Computa\c c\~ao, Universidade de S\~ao Paulo}, \orgaddress{\street{Caixa Postal 668}, \city{S\~ao Carlos}, \postcode{13560-970}, \state{SP}, \country{Brazil}}}

\affil[2]{\orgdiv{Departamento de Matem\'atica}, \orgname{Universidade Federal do Paran\'a}, \orgaddress{\street{Caixa Postal 19081}, \city{Curitiba}, \postcode{81531-990}, \state{PR}, \country{Brazil}}}


\abstract{
	We investigate the global hypoellipticity and global solvability of systems of left-invariant differential operators on compact Lie groups.  Focusing on diagonal systems, we establish necessary and sufficient conditions for these global properties. Specifically, we show that global solvability is characterized by a Diophantine condition on the symbol of the system, while global hypoellipticity further requires that a set depending on the symbol to be finite. As an application, we provide a complete characterization of these properties for systems of vector fields defined on products of tori and spheres. Additionally, we present illustrative examples, including systems involving higher-order differential operators. Finally, we extend our analysis to triangular systems on compact Lie groups, introducing an additional condition related to the boundedness of the dimensions of the representations in these Lie groups.
}



\maketitle

%============================================================
\section{Introduction}
%============================================================

We address the existence and regularity of global solutions for a class of systems of left-invariant differential operators defined on compact Lie groups. An operator or system of operators $P$ is said to be globally hypoelliptic if the smoothness of $Pu$ implies the smoothness of $u$. Similarly, $P$ is globally solvable if the equation $Pu = f$ admits a distributional solution $u$  whenever the right-hand side $f$ belongs to a suitable subset of distributions. These global properties reflect a deep interplay between the algebraic structure of the group and the analytic behavior of the operators.  

Extensive research has been conducted on these global properties for operators defined on the torus, as evidenced by works such as \cite{
	BDGK2015_jpdo,
	BP1999_jmaa,
	CH1977_pams,
	DM2016_mana,
	GW1972_pams}.
Additionally, these properties have also been studied in the context of compact Lie groups \cite{
	KMR2021_jfa,
	KKM2024,
	KMP2021_jde}, 
and some results have been obtained on compact manifolds \cite{
	ADL2023_jga,
	ADL2023_math-ann,
	AGK2018_jam,
	AK2019_jst,
	HZ2017_man,
	HZ2022_jde}.

The particular case in which the operators are vector fields defined on tori deserves special attention. As conjectured by S. Greenfield and N. Wallach (see \cite{
	Forni08_cont-math, 
	GW1973_tams}), 
if a smooth closed manifold $M$ admits a globally hypoelliptic vector field $X$, then $M$ is diffeomorphic to a torus, and $X$ can be conjugated to a constant vector field that satisfies a certain Diophantine condition.

Consequently, the investigation of global hypoellipticity for vector fields on closed manifolds has primarily focused on the torus. In this paper, we extend this investigation to diagonal systems, as defined on Section \ref{Section_Operators}. This class encompasses systems of first-order left-invariant operators, such as systems of left-invariant vector fields on products of Lie groups. Additionally, it includes classes of left-invariant operators of arbitrary order, unlike some of the previous works mentioned, which are restricted to first-order operators.

We introduce Diophantine conditions as a key tool in our analysis. These conditions provide a framework for characterizing the behavior of the symbol of the system, which is crucial for determining its global properties. The Diophantine conditions ensure a certain degree of separation between the eigenvalues of the symbol, preventing resonance phenomena that could obstruct solvability and hypoellipticity.

Our results show that the global solvability of the system is equivalent to the symbol satisfying a Diophantine condition. Furthermore, the system will be globally hypoelliptic if, in addition to a Diophantine condition, the set of singularities of the symbol is finite. These findings extend classical results in the literature and offer new perspectives for exploring this subject within the context of compact Lie groups.

This article is organized as follows. In Section 2, we provide the necessary background on Fourier analysis on compact Lie groups, establishing the foundational tools and concepts required for our results. Section 3 begins by presenting a general necessary condition for the global hypoellipticity of systems on compact Lie groups. Next, we introduce the class of diagonal systems and derive necessary and sufficient conditions for their global hypoellipticity and global solvability under various formulations. This section concludes with the study of a special subclass of diagonal systems, encompassing systems of vector fields defined on products of compact Lie groups. In Section 4, we present illustrative examples and applications that demonstrate the scope and implications of our results, focusing on diagonal systems defined on products of tori and spheres. Finally, in Section 5, we extend the analysis to triangular systems. We give examples of sufficient conditions that guarantee global hypoellipticity.

%============================================================
%============================================================
\section{Overview on Fourier Analysis in Compact Lie groups}\label{overview}
%============================================================
%============================================================

In this section, we introduce the notation and fundamental results needed for this study. A more detailed presentation of these concepts, as well as the proofs of the results discussed here, can be found in \cite{
	RT2010_book}.

Let $G$ be a compact Lie group, and let $\mu$ denote the normalized Haar measure on $G$. The set of equivalence classes of continuous irreducible unitary representations of $G$, called the unitary dual, is denoted by $\widehat{G}$. It is well known that $\widehat{G}$ is countable. By the Peter-Weyl Theorem, there exists an orthonormal basis for $L^2(G)$, which can be constructed as follows: for each class $\Xi \in \widehat{G}$, we select a representative matrix-valued function $\xi: G \to U(d_\xi)$, where $d_\xi$ is the dimension of the representation $\xi$. Writing $\xi = (\xi_{mn})_{1 \leq m,n \leq d_\xi}$, the set
$$
\bigcup_{\Xi \in \widehat{G}} \left\{ \sqrt{d_\xi} \cdot \xi_{mn}; \ 1 \leq m,n \leq d_\xi \right\}
$$
forms an orthonormal basis for $L^2(G)$. From now on, we assume that a unique representative $\xi$ has been chosen for each class $\Xi \in \widehat{G}$, although in special cases we may impose additional properties on these representatives.

Next, we can consider the Fourier analysis on $G$ with respect to this basis. For each $f \in L^1(G)$ and $\Xi = [\xi] \in \widehat{G}$, the $\xi$-Fourier coefficient of $f$ is given by the matrix
$$
\widehat{f}(\xi) \doteq \int_G f(g) \xi(g)^\ast \, d\mu(g).
$$
This is defined up to conjugation by unitary matrices, which is sufficient for our purposes, as we are primarily concerned with estimating the Hilbert-Schmidt norm of $\widehat{f}(\xi)$. 

We denote by $C^\infty(G)$ the space of smooth functions on $G$, equipped with the standard topology of uniform convergence for functions and their derivatives. The space of distributions on $G$ is denoted by $\mathcal{D}'(G)$, which is the topological dual of $C^\infty(G)$.

Let $\Delta$ be the Laplace-Beltrami operator on $G$. For each $\Xi = [\xi] \in \widehat{G}$, the matrix entries $\xi_{mn}$ are eigenfunctions of $\Delta$, all corresponding to the same eigenvalue $\lambda(\xi) \leq 0$. The operator $(I - \Delta)^{1/2}$ is positive definite, and we denote its eigenvalue corresponding to $\Xi = [\xi]$ by $\jp{\xi} \doteq (1 + \lambda(\xi))^{1/2}$.

For each $u \in \mathcal{D}'(G)$ and $\Xi = [\xi] \in \widehat{G}$, the $\xi$-Fourier coefficient of $u$ is defined by
$$
\widehat{u}(\xi) \doteq \langle u, \xi^\ast \rangle.
$$
Again, this is well-defined up to conjugation by a unitary matrix.

The Peter-Weyl basis allows us to characterize distributions and different regularity classes of functions through the behavior of their Fourier coefficients. For instance, suppose that for each $\Xi = [\xi] \in \widehat{G}$, we can associate a matrix $(x(\xi)_{mn}) \in \C^{d_\xi \times d_\xi}$, and that there exist constants $M > 0$ and $N > 0$ such that
\begin{equation}\label{dist}
	\big|x(\xi)_{mn} \big| \leq M \jp{\xi}^{N}
\end{equation}
for all $1 \leq m,n \leq d_\xi$. Then, the series
$$
u \doteq \sum_{\Xi = [\xi] \in \widehat{G}} d_\xi \sum_{m,n} x(\xi)_{mn} \xi_{nm}
$$
converges in $\mathcal{D}'(G)$ and defines a distribution $u$ such that $\widehat{u}(\xi)_{mn} = x(\xi)_{mn}$ for all $1 \leq m,n \leq d_\xi$, and $\Xi = [\xi] \in \widehat{G}$.

Conversely, if $u \in \mathcal{D}'(G)$, then there exist constants $M > 0$ and $N > 0$ such that
$$
\big|\widehat{u}(\xi)_{mn}\big| \leq M \jp{\xi}^{N}
$$
for all $1 \leq m,n \leq d_\xi$ and $\Xi = [\xi] \in \widehat{G}$.

Regarding smooth functions, a distribution $u \in \mathcal{D}'(G)$ is a smooth function if and only if, for every $N > 0$, there exists a constant $M > 0$ such that
\begin{equation}\label{smooth1}
	\big|\widehat{u}(\xi)_{mn} \big| \leq M \jp{\xi}^{-N}
\end{equation}
for all $1 \leq m,n \leq d_\xi$ and $\Xi = [\xi] \in \widehat{G}$. In this case, the corresponding Fourier series converges in the $L^2(G)$-norm, and the Plancherel formula holds.

For a continuous operator $P: \mathcal{D}'(G) \to \mathcal{D}'(G)$, the symbol of $P$ is defined at $g \in G$ and $\Xi = [\xi] \in \widehat{G}$ by the matrix
$$
\sigma_P(g,\xi) \doteq \xi(g)^\ast (P\xi)(g),
$$
where $P\xi$ is the matrix $(P\xi)_{mn} \doteq P(\xi_{mn})$ for $1 \leq m,n \leq d_\xi$. This is also defined up to conjugation by a unitary matrix. In the special case where $P$ is left-invariant (i.e., it commutes with left translations), $\sigma_P$ does not depend on $g \in G$, and for all $u \in \mathcal{D}'(G)$ and $\Xi = [\xi] \in \widehat{G}$, the following formula holds:

\begin{equation} \label{multiply}
\widehat{Pu}(\xi) = \sigma_P(\xi)\cdot \widehat{u}(\xi).
\end{equation}

%============================================================
%============================================================	
\section{Systems of operators}\label{Section_Operators}
%============================================================
%============================================================

The first result of this paper establishes a general property that applies to any system of left-invariant continuous operators on a compact Lie group. This result provides a sufficient condition for global hypoellipticity, which will be used in subsequent results. Next, we introduce the notion of diagonal systems and analyze their properties.

Let $ P_1, \dots, P_r: \mathcal{D}'(G) \to \mathcal{D}'(G) $ be left-invariant continuous operators, and consider the system $ P \doteq (P_1, \dots, P_r): \mathcal{D}'(G) \to \mathcal{D}'(G)^r $. The total symbol of $ P $ is defined as
$$
\sigma_P(\xi) \doteq (\sigma_{P_1}(\xi), \dots, \sigma_{P_r}(\xi)),
$$
and its norm is given by
$$
\|\sigma_P(\xi)\| \doteq \max_{1 \leq j \leq r} \left\| \sigma_{P_j} (\xi) \right\|_{0},
$$
where $ \left\|A \right\|_{0} = \max_{m,n} |a_{mn}| $ for any matrix $ A = (a_{mn}) \in \C^{d_\xi \times d_\xi} $.

Notice that if $P$ is a such a system and $u \in \mathcal{D}'(G)$ and $f=(f_1,...,f_r) \in \mathcal{D}'(G)^r$ are such that $Pu=f$, that is, $P_ju=f_j$ for $j=1,...,r$, then for each $[\xi] \in \widehat{G}$ we have $\widehat{P_j u}(\xi) = \widehat{f_j}(\xi)$. By property (\ref{multiply}), we have $\sigma_{P_j}(\xi) \cdot \widehat{u}(\xi) = \widehat{f_j}(\xi)$ for all $j=1,...,r$. In particular, the vector of matrices $\widehat{f}(\xi):=(\widehat{f_1}(\xi),...,\widehat{f_r}(\xi))$ lies in $Im(\sigma_P(\xi))$, namely, the image of the operator $\sigma_P(\xi)$. In this way, we will define the space

\begin{equation}\label{A(P)}
	\mathcal{A}(P) \doteq \{ f\in \mathcal{D}'(G)^r; \forall [\xi]\in \widehat{G}, \widehat{f}(\xi) \in Im(\sigma_P(\xi)) \}.
\end{equation}

We will call $\mathcal{A}(P)$ the space of compatibility conditions. Later we will consider this space in order to define our notion of global solvability. But before, we will study regularity of the system $P$ in the following sense:

\begin{definition}
	The system $P$ is said to be globally hypoelliptic (or simply (GH)) if the conditions 
	$u \in \mathcal{D}'(G)$ and $Pu \in C^\infty(G)^r$ imply that $u \in C^\infty(G).$
\end{definition}

\begin{theorem}\label{Zinfinite-notGH}
	If $P$ is (GH) then the following set is finite:
	\begin{equation}\label{zset} 
		\mathcal{Z} \doteq \left\{ \Xi = [\xi] \in \widehat{G}; \ \bigcap_{j=1}^r \mathrm{Ker}(\sigma_{P_j}(\xi)) \neq 0 \right\}.
	\end{equation}
\end{theorem}
\begin{proof}
	Suppose that $\mathcal{Z}$ is infinite. In this case, for each $k \in \mathbb{N}$, there exists $\Xi_k = [\xi_k] \in \widehat{G}$ such that 
	$$
	\bigcap_{j=1}^r \mathrm{Ker}(\sigma_{P_j}(\xi_k)) \neq 0.
	$$
	
	For each $k \in \mathbb{N}$, we can select a matrix $x(\xi_k) \in \C^{d_{\xi_k} \times d_{\xi_k}}$ such that
	$$
	\left\|x(\xi_k)\right\|_{0} = 1 \quad \text{and} \quad P_j x(\xi_k) = 0 \quad \text{for all } j = 1, \dots, r.
	$$
	
	Since the sequence of matrices $(x(\xi_k))_{k \in \mathbb{N}}$ satisfies condition (\ref{dist}), it defines a distribution $u \in \mathcal{D}'(G)$. Since 	$\|x(\xi_k)\|_{0} = 1$ for all $k$, this distribution is not induced by a smooth function on $G$. Morever, by construction, $Pu = 0 \in C^\infty(G)$. This implies that $P$ is not (GH).
\end{proof} 

%====================================================
%====================================================
\subsection{Global hypoellipticity for diagonal systems} 
%====================================================
%====================================================

In this subsection, we provide a complete characterization of global hypoellipticity for diagonal systems of left-invariant operators on compact Lie groups.

\begin{definition}
	A system of left-invariant continuous operators $P = (P_1, \dots, P_r)$ is said to be diagonal if, for each $\Xi \in \widehat{G}$, there exists a matrix representative $\xi \in \Xi$ such that the matrix $\sigma_{P_j}(\xi)$ is diagonal for all $1 \leq j \leq r$. 	
	In this case, we write
	$$
	\sigma_{P_j}(\xi) = \mathrm{diag}(\sigma_1^j(\xi), \dots, \sigma_{d_\xi}^j(\xi))
	$$
	for all $\Xi = [\xi] \in \widehat{G}$ and $1 \leq j \leq r$.	
\end{definition}

Note that for a such diagonal system $P$, the set $\mathcal{Z}$ defined in (\ref{zset}) takes the form:
$$
\mathcal{Z} = \{\Xi=[\xi] \in \widehat{G}; \exists k_0 \in \{1,...,d_\xi\} \textrm{ such that } \sigma_{k_0}^j(\xi)=0, \forall j=1,...,r \}.
$$

We will present examples of diagonal systems in the section \ref{examples}. In the context of constant coefficient operators, it is usual to ``Diophantine condition'' appear. They usually tell how fast certain quantities can be approximated by polinomials, exponentials or other kinds of growths. In this work, we will define the following condition.

\begin{definition}\label{DC-diagonal-def}
	A diagonal system $P = (P_1, \dots, P_r)$ is said to satisfy the condition \eqref{DC1} if there exist constants $C > 0$ and $M \in \mathbb{R}$ such that
	\begin{equation}\label{DC1}
		\max_{j \in \bigwedge_{k,\xi}} \big|\sigma_k^j(\xi) \big| \geq C \cdot \jp{\xi}^{-M}. \tag{DC}
	\end{equation}
	whenever $\bigwedge_{k,\xi}:=\{j \in \{1,...,r\}; \sigma_k^j(\xi) \neq 0\} \neq \emptyset$.
\end{definition}

The first result characterize the global hypoellipticity of a diagonal system in terms of the condition (\ref{DC1}) and the set defined in (\ref{zset}).

\begin{theorem}\label{ghthm}
	A diagonal system $P = (P_1, \dots, P_r)$ is (GH) if and only if the set $\mathcal{Z}$ is finite and the condition \eqref{DC1} holds.
\end{theorem}
\begin{proof}
	Assume first that $\mathcal{Z}$ is finite and  \eqref{DC1} holds. Given $f = (f_1, \dots, f_r) \in \mathcal{A}(P) \cap C^\infty(G)^r$, suppose there exists $u \in \mathcal{D}'(G)$ such that $Pu = f$. We will prove that $u \in C^\infty(G)$.
	
	Let $C > 0$ and $M  \in \mathbb{R}$ be the constants such that condition (\ref{DC1}) holds and fix $N>0$. Since we are assuming that $\mathcal{Z}$ is finite, tha quantity
	$$
	C_1(N) \doteq \max \left\{ \big|\widehat{u}(\xi)_{k\ell} \big|\jp{\xi}^{N}; \ [\xi] \in \mathcal{Z}, \ 1 \leq k,\ell \leq d_\xi \right\}
	$$
	is a non-negative real number. By definition of $C_1(N)$, we have
	$$\big|\widehat{u}(\xi)_{k\ell} \big| \leq C_1(N) \jp{\xi}^{-N}$$
	for all $[\xi]\in \mathcal{Z}$.
	
	Now fix $[\xi] \in \widehat{G}\setminus \mathcal{Z}$ and $k_0,l_0$ with $1 \leq k_0,l_0 \leq d_\xi$. We want to estimate $|\widehat{u}(\xi)_{k_0 l_0}|$. Since $[\xi] \in \widehat{G}\setminus \mathcal{Z}$, we have $\max_j |\sigma_{k_0}^j(\xi)|>0$, so by the condition (\ref{DC1}) there exists $j_0$ with
	$$\big|\sigma_{k_0}^{j_0}(\xi) \big| \geq C \cdot \jp{\xi}^{-M}.$$

	Since $P_{j_0} u = f_{j_0}$ and $\sigma_{P_{j_0}}(\xi)$ is diagonal, we have
	$$
	\widehat{u}(\xi)_{k_0 \ell_0} =  (\sigma_{k_0}^{j_0}(\xi))^{-1}\widehat{f}_{j_0}(\xi)_{k_0 \ell_0}.
	$$
	
	Given that $f_{j_0} \in C^\infty(G)$, there exists $C_2(N) > 0$ such that
	$$
	\big|\widehat{f}_{j_0}(\xi)_{k \ell}\big| \leq C_2(N) \jp{\xi}^{-N - M}
	$$
	for all $[\xi] \in \widehat{G}$ and $1 \leq k, \ell \leq d_\xi$.
	
	Thus, we have
	$$
	\big|\widehat{u}(\xi)_{k_0 \ell_0} \big| \leq \frac{1}{C} \jp{\xi}^{M} \cdot C_2 \jp{\xi}^{-N - M} = \frac{C_2(N)}{C} \cdot \jp{\xi}^{-N}.
	$$
	
	Letting $\widetilde{C} \doteq \max \{C_1(N), C_2(N) / C\} > 0$, then $\widetilde{C}$ depends only on $N$ and we obtain
	$$
	\big|\widehat{u}(\xi)_{k \ell} \big| \leq \widetilde{C} \cdot \jp{\xi}^{-N}
	$$
	for all $[\xi] \in \widehat{G}$ and $1 \leq k, \ell \leq d_\xi$.
	This shows that $u \in C^\infty(G)$, proving that $P$ is (GH). 
	
	Conversely, suppose that the condition \eqref{DC1} does not hold. In this case, for each $n \in \mathbb{N}$, there exist $\Xi_n = [\xi_n] \in \widehat{G}$ and $1 \leq k_n \leq d_{\xi_n}$, with $\Xi_n \neq \Xi_m$ for all $n \neq m$, such that $\max_j |\sigma_{k_n}^{j}(\xi_n)|>0$ and for all $j=1,...,r$ the following inequality holds:
	$$
	\big|\sigma_{k_n}^{j}(\xi_n)\big| < \jp{\xi_n}^{-n}.
	$$
	For each $j = 1, \dots, r$, we define a sequence of matrices $(x^j(\xi))_{[\xi] \in \widehat{G}}$ as follows:
	$$
	x^j(\xi)_{k\ell} \doteq \begin{cases}
		\sigma_{k_n}^{j}(\xi_n), & \text{if } \Xi = \Xi_n, \text{ and } (k,\ell) = (k_n, k_n), \\
		0, & \text{otherwise}.
	\end{cases}
	$$
	
	In this way, for each $j = 1, \dots, r$, the sequence $(x^j(\xi))_{[\xi] \in \widehat{G}}$ defines a function $f_j \in C^\infty(G)$. To see this, first fix $j \in \{1,...,r\}$ and $N > 0$. By definition of $x^j(\xi)$, it is sufficient to estimate $|x^j(\xi_n)_{k_n k_n}|$ for $n \in \mathbb{N}$. So fix $n \in \mathbb{N}$. If $ n\leq N$, we define 
	$$C_1(N):= \{|x^j(\xi_m)_{k_m k_m}|\jp{\xi_m}^N; m \leq N \} \geq 0$$ 
	and then we get $|x^j(\xi_n)_{k_n k_n}| \leq C_1(N) \jp{\xi_n}^{-N}$.
	
	If $n\geq N$, then
	$$\big|x^{j}(\xi_n)_{k_n k_n} \big| = \big|\sigma_{k_n}^{j} (\xi_n)\big| < \jp{\xi_n}^{-n} < \jp{\xi_n}^{-N}.
	$$
	
	If we define $\widetilde{C}(N) \doteq \max \{C_1(N), 1\} > 0$, then for all $[\xi] \in \widehat{G}$ and $1 \leq k,\ell \leq d_{\xi}$, we have
	$$
	\big|x^j(\xi)_{k\ell} \big| \leq \widetilde{C}(N) \jp{\xi}^{-N},
	$$
	which guarantees that $(x^j(\xi))_{[\xi] \in \widehat{G}}$ defines a smooth function $f_j$ for each $1 \leq j \leq r$. It is easy to see that $f = (f_1, \dots, f_r) \in \mathcal{A}(P)$. Suppose that there exists $u \in \mathcal{D}'(G)$ such that $Pu = f$. Since $\max_j |\sigma_{k_n}^j(\xi_n)|>0$, there exists $j_n \in \{1,...,r\}$ such that
	
	$$
	0<\big|\sigma_{k_n}^{j}(\xi_n)\big| < \jp{\xi_n}^{-n}.
	$$
	Then, for each $n \in \mathbb{N}$, we have $P_{j_n} u = f_{j_n}$ and so
	$$
	\sigma_{k_n}^{j_n}(\xi_n) \cdot \widehat{u}(\xi_n)_{k_n k_n} = \widehat{f_{j_n}}(\xi_n)_{k_n k_n} = \sigma_{k_n}^{j_n}(\xi_n).
	$$
	Since $\sigma_{k_n}^{j_n}(\xi_n) \neq 0$, we have $\widehat{u}(\xi_n)_{k_n k_n} = 1$ for all $n \in \mathbb{N}$.
	
	Since the set $\{[\xi_n]; n \in \mathbb{N}\}$ is infinite, such $u$ cannot be a smooth function, and therefore $P$ is not (GH). This concludes the proof of the theorem.
\end{proof} 

If $P = (P_1, \ldots, P_r)$ is a system of operators on $G$, it is clear that if any of the operators $P_j$ is (GH), then the entire system $P$ is also (GH). Indeed, suppose that $f = (f_1, \ldots, f_r) \in C^\infty(G)^r$ and $u \in \mathcal{D}'(G)$ are such that $Pu = f$. In particular, this implies that $P_ju = f_j$. Since $f_j \in C^\infty(G)$ and $P_j$ is (GH), it follows that $u \in C^\infty(G)$. However, the converse is not true, as we illustrate in the following example.

\begin{example}
	Let $G$ be a compact Lie group, and suppose that $\widehat{G} = \{[\xi_n] : n \in \mathbb{N}\}$ with $d_{\xi_n} = 1$ for all $n \in \mathbb{N}$. This is the case when $G$ is any compact abelian Lie group. Now, let $P_1, P_2$ be left-invariant operators on $G$ such that their symbols satisfy
	$$
	\sigma_{P_1}(\xi_n) = 
	\begin{cases} 
		\jp{\xi_n} & \text{if } n \text{ is even}, \\ 
		0 & \text{if } n \text{ is odd},
	\end{cases}
	$$
	and
	$$
	\sigma_{P_2}(\xi_n) = 
	\begin{cases} 
		0 & \text{if } n \text{ is even}, \\ 
		\jp{\xi_n} & \text{if } n \text{ is odd}.
	\end{cases}
	$$
	Consider the corresponding system $P = (P_1, P_2)$. In this case, the system is trivially diagonal. By construction, it is clear that $\mathcal{Z} = \emptyset$ and condition (\ref{DC1}) is satisfied. Therefore, by Theorem \ref{ghthm}, the system $P$ is (GH). Nevertheless, for each $j = 1, 2$, the symbol $\sigma_{P_j}(\xi_n)$ vanishes for infinitely many values of $n$, implying that $P_j$ is not (GH).
\end{example}

%====================================================
%====================================================
\subsection{Global solvability for diagonal systems} 
%====================================================
%====================================================

In this subsection, we provide a comprehensive characterization of global solvability for diagonal systems of left-invariant operators, addressing both distributional and smooth solvability.

Let $P = (P_1, \ldots, P_r)$ be a diagonal system on $G$, and assume that $u \in \mathcal{D}'(G)$ and $f=(f_1, f_2, \dots, f_r) \in \mathcal{D}'(G)^r$ satisfy the equation $Pu = f$. Since the symbols $\sigma_{P_1}(\xi), \dots, \sigma_{P_r}(\xi)$ are diagonal, they commute. Therefore, we have 
$$
P_j f_\ell = P_j P_\ell u = P_\ell P_j u = P_\ell f_j
$$
for all $1 \leq j, \ell \leq r$. As a consequence, we obtain the following relation:
\begin{equation}\label{comp1}
	\sigma_{P_j}(\xi) \cdot \widehat{f_\ell}(\xi) = \sigma_{P_\ell}(\xi) \cdot \widehat{f_j}(\xi)
\end{equation}
for all $\Xi = [\xi] \in \widehat{G}$ and $1 \leq j, \ell \leq r$.

Moreover, since $\sigma_{P_j}(\xi) \cdot \widehat{u}(\xi) = \widehat{f_j}(\xi)$, we have $\sigma_k^j(\xi) \cdot \widehat{u}(\xi)_{k\ell} = \widehat{f_j}(\xi)_{k\ell}$,  for all $1 \leq k, \ell \leq d_\xi$ and $\Xi = [\xi] \in \widehat{G}$. Thus,
\begin{equation}\label{comp2}
	\sigma_k^j(\xi) = 0 \Rightarrow \widehat{f_j}(\xi)_{k\ell} = 0,  \quad 1 \leq \ell \leq d_\xi.
\end{equation}

Next, we define our notion of global solvability.

\begin{definition}
	We say that $P$ is globally solvable (or just (GS)) if 
	$$
	P( \mathcal{D}'(G)) = \mathcal{A}(P).
	$$
\end{definition}

\begin{theorem}\label{gs_condition}
	A diagonal system $P$ is (GS) if, and only if, $P$ satisfies the condition \eqref{DC1}.
\end{theorem}

\begin{proof}
	First, suppose that $P$ does not satisfy the condition \eqref{DC1}. In this case, for each $n \in \mathbb{N}$, there exist $\Xi_n = [\xi_n] \in \widehat{G}$ and $1 \leq k_n \leq d_{\xi_n}$, with $\Xi_n \neq \Xi_m$ for all $n \neq m$ and $\max_j |\sigma_{k_n}^j(\xi_n)|>0$, such that
	$$
	\big|\sigma_{k_n}^{j}(\xi_n) \big| < \jp{\xi_n}^{-n}
	$$
	for all $n \in \mathbb{N}$. For each fixed $j \in \{1,...,r\}$, we define a sequence of matrices $(x^j(\Xi))_{[\xi] \in \widehat{G}}$ by
	$$
	x^j(\xi)_{k\ell} \doteq \begin{cases} 
		\sigma_{k_n}^j(\xi_n) \jp{\xi_n}^n, & \text{if } [\xi] = [\xi_n] \text{ and } (k, \ell) = (k_n, k_n), \\
		0, & \text{otherwise}. 
	\end{cases}
	$$
	We claim that the sequence $(x^j(\xi))_{[\xi] \in \widehat{G}}$ defines a distribution $f_j \in \mathcal{D}'(G)$. In fact, as before, it is sufficient to estimate $|x^j(\xi_n)_{k_n k_n}|$. But in this case we have
	$$|x^j(\xi_n)_{k_n k_n}| = |\sigma_{k_n}^j(\xi_n)| \jp{\xi_n}^n \leq \jp{\xi_n}^{-n}\cdot \jp{\xi_n}^n =1=\jp{\xi_n}^0$$
	for all $n \in \mathbb{N}$. It is also clear that $f:=(f_1,...,f_r) \in \mathcal{A}(P)$. We claim that there is no $u \in \mathcal{D}'(G)$ such that $Pu = f$. Indeed, suppose, by contradiction, that such a $u$ exists. Since $\max_j |\sigma_{k_n}^j(\xi_n)|>0$, there exists $j_n \in \{1,...,r\}$ such that 
	$$
	0<\big|\sigma_{k_n}^{j}(\xi_n) \big| < \jp{\xi_n}^{-n}
	$$
	for all $n \in \mathbb{N}$. Since $Pu=f$, we must have $P_{j_n} u=f_{j_n}$, and since $\sigma_{P_{j_n}}(\xi_n)$ is diagonal, we get
	$$\sigma_{k_n}^{j_n}(\xi_n)\cdot \widehat{u}(\xi_n)_{k_n k_n}  = \widehat{f}_{j_n}(\xi_n)_{k_n k_n} \Rightarrow \widehat{u}(\xi_n)_{k_n k_n}=\jp{\xi_n}^n$$
	for all $n \in \mathbb{N}$. In this way, the sequence $(\widehat{u}(\xi))_{[\xi]\in \widehat{G}}$ cannot define a distribution. This concludes that $P$ is not (GS).
	
	Conversely, let $C_1>0$ and $N_1 \in \mathbb{R}$ be constants such that (\ref{DC1}) holds. Given $f = (f_1, \dots, f_r) \in \mathcal{A}(P)$, we need to find a distribution $u \in \mathcal{D}'(G)$ such that $Pu = f$. 
	
	Fix $\Xi = [\xi] \in \widehat{G}$ and $1 \leq k, \ell \leq d_\xi$. If $\sigma_k^j(\xi) = 0$ for all $j\in \{1,..,r\}$, we define $\widehat{u}(\xi)_{k\ell}\doteq 0$. If there exists $j_0 \in \{1,...,r\}$ such that $\sigma_k^{j_0}(\xi)\neq 0$, we define
	$$
	\widehat{u}(\xi)_{k\ell} \doteq \frac{\widehat{f_{j_0}}(\xi)_{k\ell}}{\sigma_k^{j_0}(\xi)}.
	$$
	Since $f \in \mathcal{A}(P)\subset \widetilde{\mathcal{A}}(P)$, the relations (\ref{comp1}) guarantees that the number above is well-defined and does not depend on $j_0$. The identity $\widehat{(P_j u)}(\xi)_{k\ell} = (\sigma_{P_j} \cdot \widehat{u}(\xi))_{k\ell} = \widehat{f_j}(\xi)_{k\ell}$ holds by construction for all $j = 1, \dots, r$. 
	
	Given that $f_{j_0} \in \mathcal{D}'(G)$, there are constants $C_2>0$ and $N_2\in \mathbb{R}$ such that
	$$
	\big|\widehat{f_{j_0}}(\xi)_{k\ell} \big| \leq C_2 \jp{\xi}^{N_2}
	$$
	for all $\Xi = [\xi] \in \widehat{G}$ and $1 \leq k, \ell \leq d_\xi$. Thus,
	$$
	\big|\widehat{u}(\xi)_{k\ell} \big| \leq C_1 C_2 \jp{\xi}^{N_1 + N_2},
	$$
	and the sequence $(\widehat{u}(\xi))_{[\xi] \in \widehat{G}}$ indeed defines a distribution $u \in \mathcal{D}'(G)$, which, by construction, satisfies $Pu = f$. Therefore, $P$ is (GS), and this concludes the proof.
\end{proof}

An immediate corollary is the following.

\begin{corollary}
	If $P$ is a (GH) diagonal system, then it is also (GS).
\end{corollary}

Finally, observe that the results discussed above provide a solution to the equation $Pu = f$ only in the distributional sense. We now turn our attention to the smooth solvability.

\begin{definition}
	A diagonal system $P$ is said to be $C^\infty$-globally solvable (or simply $C^\infty$-(GS)) if
	$$
	P(C^\infty(G)) = \mathcal{A}(P) \cap C^\infty(G)^r.
	$$
\end{definition}

Assume that the condition \eqref{DC1} holds, and let $$f = (f_1, \dots, f_r) \in \mathcal{A}(P) \cap C^\infty(G)^r.$$ 

Let $C_1 > 0$ and $N_1 \in \mathbb{R}$ be the constants defined in condition \eqref{DC1}. From the proof of the previous theorem, we know that there exists a distribution $u \in \mathcal{D}'(G)$ such that $Pu = f$. We now claim that, in this special case, the distribution $u$ constructed in the previous proof is in fact a smooth function.

To see this, fix $N > 0$, $\Xi_0 = [\xi_0] \in \widehat{G}$, and $1 \leq k_0, \ell_0 \leq d_{\xi_0}$. If $\sigma_{k_0}^j(\xi_0) = 0$ for all $j \in \{1, \dots, r\}$, then by definition, $\widehat{u}(\xi_0)_{k_0 \ell_0} = 0$, and no further estimate is needed in this case.

Now, suppose that $\sigma_{k_0}^{j_0}(\xi_0) \neq 0$ for some $j_0 \in \{1, \dots, r\}$. Since $f_{j_0} \in C^\infty(G)$, there exists a constant $C_{j_0, 2}$, depending only on $N$ and $j_0$, such that
$$
\big|\widehat{f_{j_0}}(\xi)_{k \ell}\big| \leq C_{j_0, 2} \langle \xi \rangle^{-N_1 - N}
$$
for all $\Xi = [\xi] \in \widehat{G}$ and $1 \leq k, \ell \leq d_\xi$. Thus, we obtain
$$
\big|\widehat{u}(\xi)_{k_0 \ell_0}\big| = \frac{1}{|\sigma_{k_0}^{j_0}(\xi)|} \cdot \big|\widehat{f_{j_0}}(\xi)_{k_0 \ell_0}\big| \leq \frac{1}{C_1} \cdot \langle \xi \rangle^{N_1} \cdot C_{j_0, 2} \langle \xi \rangle^{-N_1 - N}.
$$

If we define $C = \max_j \left\{C_j, \frac{2}{C_1}\right\} > 0$, then $C$ is a positive constant depending only on $N$, and we have
$$
\big|\widehat{u}(\xi_0)_{k_0 \ell_0}\big| \leq C \langle \xi_0 \rangle^{-N}.
$$
This shows that $u \in C^\infty(G)$, leading to the following corollary.

\begin{corollary}\label{corollary_AGH}
	If $P$ is a diagonal system and the condition \eqref{DC1} holds, then $P$ is $C^\infty$-globally solvable.
\end{corollary}

\begin{remark}\label{AGH->GS}
	Assuming that condition (\ref{DC1}) is satisfied, the preceding arguments leading to Corollary \ref{corollary_AGH} reveal the following property: for any distribution $v \in \mathcal{D}'(G)$ such that $Pv \in C^\infty(G)$, there exists a smooth function $u \in C^\infty(G)$ satisfying $Pu = Pv$. 
	The operators with this property are known as \textit{almost-globally hypoelliptic}. Hence, it follows from \cite[Theorem 2.2]{AFJR2024_math-ann}, if our operator $P$ satisfies (DC) then $P$ has closed image, which coincides with our notion of global solvability.
\end{remark}

%====================================================
%====================================================
\subsection{Systems of vector fields on product of compact Lie groups} 
%====================================================
%====================================================

In this subsection, as a consequence of the results in the previous subsections, we present a complete characterization of global hypoellipticity and global solvability for this important class of diagonal systems, consisting of left-invariant vector fields defined on the Cartesian product of compact Lie groups.

Let $G_1, \dots, G_n$ be compact Lie groups, and consider the product Lie group $\mathcal{G} \doteq G_1 \times \dots \times G_n$. It is well known that $\widehat{\mathcal{G}} \simeq \widehat{G_1} \times \dots \times \widehat{G_n}$, meaning that for each class $\Xi \in \widehat{\mathcal{G}}$, there exist unitary, irreducible, and continuous representations $\xi^1, \dots, \xi^n$ of $G_1, \dots, G_n$, respectively, such that $\Xi = [\xi^1 \otimes \dots \otimes \xi^n]$.

In view of the Peter-Weyl theorem, for each $\Xi \in \widehat{\mathcal{G}}$ and $j = 1, \dots, n$, we fix a matrix representative $\xi^j$ of $G_j$ such that $\Xi = [\xi^1 \otimes \dots \otimes \xi^n]$. And, for each $\Xi \in \widehat{\mathcal{G}}$, define the set 
$$
J_\Xi \doteq \left\{(a_1, \dots, a_n) \in \mathbb{N}^n; 1 \leq a_j \leq d_{\xi^j} \text{ for } j = 1, \dots, n \right\}.
$$ 

Thus, for each $u \in \mathcal{D}'(\mathcal{G})$, $\Xi \in \widehat{\mathcal{G}}$, and $A, B \in J_\Xi$, the corresponding Fourier coefficient is given by
$$
\widehat{u}(\Xi)_{AB} \doteq \jp{u, \overline{\xi_{BA}}},
$$
where $\xi_{BA} \doteq \xi_{b_1 a_1}^1 \cdots \xi_{b_n a_n}^n$. The distribution $u$ can then be recovered as the series
$$
u = \sum_{\Xi \in \widehat{\mathcal{G}}} d_\Xi \sum_{A, B \in J_\Xi} \widehat{u}(\Xi)_{AB} \cdot \xi_{BA}.
$$

With this notation, the characterizations of distributions and smooth functions presented in Section \ref{overview} translate into estimating the coefficients $\widehat{u}(\Xi)_{AB}$ for $A,B \in J_{\Xi}$, while replacing $\jp{\xi}$ by $\jp{\xi^1} + \dots + \jp{\xi^n}$.

Next, let $X_j$ be a left-invariant vector field on $G_j$, for $j = 1, \dots, n$. Recall that the operators $i X_j$ are symmetric when acting on the space $L^2(G_j)$ with respect to the standard $L^2$ inner product, which implies that the symbol $\sigma_{iX_j}$ is diagonalizable.

From now on, for each $\Xi \in \widehat{\mathcal{G}}$, we fix a representative $[\xi^j] \in \widehat{G_j}$ such that $\Xi = [\xi] = [\xi^1 \otimes \dots \otimes \xi^n]$, and
$$
\sigma_{X_j}(\xi) = \mathrm{diag}(i\mu_1^j(\xi), \dots, i\mu_{d_{\xi^j}}^j(\xi)),
$$
where $\mu_k^j(\xi)$ are real numbers for all $j$ and $k$.

\begin{remark}
	As mentioned in the introduction, it was conjectured by S. Greenfield and N. Wallach in 1973 (see \cite{
		Forni08_cont-math,
		GW1973_tams}) that if a smooth closed manifold admits a globally hypoelliptic vector field, then this manifold is diffeomorphic to a torus, and the vector field can be conjugated to a constant vector field satisfying a Diophantine condition. Consequently, the investigation of global hypoellipticity for vector fields on closed manifolds has primarily focused on the torus. Notable references in this context include \cite{
		CC2000_cpde, 
		DGY2002_pems, 
		GW1972_pams}, 
	among others.
	
	In light of the likely validity of this conjecture, and to extend the scope of analysis beyond systems of vector fields on tori, a natural approach is to consider lower-order perturbations of these operators, as was done in \cite{
		KMR2020_bsm, 
		KMR2021_jfa, 
		KKM2024, 
		KMP2021_jde} 
	and other references.
\end{remark}

For each $j = 1, \dots, r$, define an operator $P_j : \mathcal{D}'(\mathcal{G}) \to \mathcal{D}'(\mathcal{G})$ by
$$
P_j \doteq c_{j1}X_1 + \dots + c_{jn}X_n + q_j,
$$
where $(c_{k\ell}) \in \mathbb{C}^{r \times n}$ and $q_1, \dots, q_r \in \mathbb{C}$. We also consider the system of operators $P = (P_1, \dots, P_r) : \mathcal{D}'(\mathcal{G}) \to \mathcal{D}'(\mathcal{G})^r$.

Note that the vector fields $X_1, \dots, X_n$ are defined on different factors of $\mathcal{G}$, and therefore, $[P_j, P_k] = 0$ for all $j, k = 1, \dots, n$.

Suppose that $f \in \mathcal{D}'(\mathcal{G})^r$ and $u \in \mathcal{D}'(\mathcal{G})$ are such that $Pu = f$. By using basic properties of Fourier analysis and the notation introduced above, we have
$$
i(c_{j1} \mu_{a_1}^1(\xi) + \dots + c_{jn} \mu_{a_n}^n(\xi) - iq_j) \widehat{u}(\xi)_{AB} = \widehat{f_j}(\xi)_{AB}
$$
for all $\Xi \in \widehat{\mathcal{G}}$, $A, B \in J_\Xi$, and $j \in \{1, \dots, r\}$, where $A = (a_1, \dots, a_n)$.

By writing $c_j \doteq (c_{j1}, \dots, c_{jn})$ and $\mu_A(\xi) \doteq (\mu_{a_1}^1(\xi), \dots, \mu_{a_n}^n(\xi))$, we can rewrite the equations above as
$$
i \left(\jp{c_j, \mu_A(\xi)} - iq_j \right) \cdot \widehat{u}(\xi)_{AB} = \widehat{f_j}(\xi)_{AB}
$$
for all $\Xi \in \widehat{\mathcal{G}}$, $A, B \in J_\Xi$, and $j \in \{1, \dots, r\}$.

Next, we translate the condition (\ref{DC1}) in this case. The system $P$ satisfies the condition (\ref{DC1}) if there are constants $C>0$ and $N \in \mathbb{R}$ with the following property:
for all $[\xi] \in \widehat{\mathcal{G}}$ such that there exists $A \in J_{[\xi]}$ with $\max_j |\jp{c_j, \mu_A(\xi)} - iq_j| > 0$, then there exists $j \in \{1,...,r\}$ such that
\begin{equation}\label{DC2}
	\big|\jp{c_j, \mu_A(\xi)} - iq_j \big| \geq C \cdot \left(\jp{\xi^1} + \dots + \jp{\xi^n} \right)^{-N} \tag{DC$_\mathcal{G}$}.
\end{equation}

Finally, the set $\mathcal{Z}$ defined in (\ref{zset}) takes the following form in this context:
\begin{align}
	\mathcal{Z}_\bullet=\{ \Xi=[\xi] \in \widehat{G}; \exists A \in J_\Xi \textrm{ such that } \jp{c_j,\mu_A(\xi)}-iq_j=0, \forall j=1,...,r\}.
\end{align}

With the notation, Theorems \ref{ghthm} and \ref{gs_condition} have the following version:

\begin{proposition}\label{ghgs-thm}
	Let $G_1, \dots, G_n$ be compact Lie groups, and $\mathcal{G} \doteq G_1 \times \dots \times G_n$.
	Consider the system $P = (P_1, \dots, P_r) : \mathcal{D}'(\mathcal{G}) \to \mathcal{D}'(\mathcal{G})^r$, where 
	$$
	\begin{cases}
		P_1 =  c_{11}X_1 + \dots + c_{1n}X_n + q_1, \\
		P_2 =  c_{21}X_1 + \dots + c_{2n}X_n + q_2,\\
		\hspace{6mm} \vdots   \\
		P_r =  c_{r1}X_1 + \dots + c_{rn}X_n + q_r,
	\end{cases}
	$$
	with $(c_{k\ell}) \in \mathbb{C}^{r \times n}$, $q_1, \dots, q_r \in \mathbb{C}$, and where each $X_j$ is a left-invariant vector field defined on $G_j$, for $j = 1, \dots, n$. Then:
	\begin{enumerate}
		\item $P$ is (GH) if, and only if, $P$ satisfies the condition \eqref{DC2} and $\mathcal{Z}_\bullet$ is finite;
		\item $P$ is (GS) if, and only if, $P$ satisfies the condition \eqref{DC2}.
	\end{enumerate}
\end{proposition}

%============================================================
%============================================================	
\section{Application: Diagonal Systems on Products of Tori and Spheres}\label{examples}
%============================================================
%============================================================	

A very interesting non-commutative Lie group, whose unitary dual can be explicitly described, is the 3-sphere $\mathbb{S}^3 \doteq \{x \in \mathbb{R}^4; \|x\|_2 = 1\} \equiv SU(2)$. A classic result shows that the unitary dual $\widehat{\mathbb{S}^3}$ is in bijection with the set $\frac{1}{2}\mathbb{N}_0$, where for each $\ell \in \frac{1}{2}\mathbb{N}_0$, there is a unique, up to isomorphism, continuous irreducible unitary representation $t^\ell : \mathbb{S}^3 \to U(2\ell + 1)$. It is common to represent the entries of the matrix-valued function $t^\ell$ by $t^\ell_{mn}$, where $m, n \in J_\ell \doteq \{-\ell, -\ell + 1, -\ell + 2, \dots, \ell - 1, \ell\}$.

The Lie algebra $\mathfrak{su}(2)$ has a standard basis $\{Y_1, Y_2, Y_3\}$, which satisfies the commutation relations $[Y_1, Y_2] = Y_3$, $[Y_2, Y_3] = Y_1$, and $[Y_3, Y_1] = Y_2$. The left-invariant operators associated with these Lie algebra elements will be denoted by $D_1$, $D_2$, and $D_3$, respectively. It can be shown that $D_3$ satisfies
$$
D_3(t^\ell_{mn}) = -in \cdot t^\ell_{mn}
$$
for all $m, n \in J_\ell$ and $\ell \in \frac{1}{2}\mathbb{N}_0$. 

The action of the operators $D_1$ and $D_2$ on the functions $t^\ell_{mn}$ can also be expressed, but the formulas are more complicated. However, there is an alternative basis for $\mathfrak{su}(2)$ (over $\mathbb{C}$), denoted by $\{\partial_0, \partial_+, \partial_-\}$, where their action on the functions $t^\ell_{mn}$ becomes much simpler. Defining
$$
\partial_+ \doteq iD_1 - D_2, \quad \partial_- \doteq iD_1 + D_2, \quad \partial_0 \doteq iD_3,
$$
it is well known that:
\begin{align}\label{tabela}
	\partial_+(t^\ell_{mn})  = & -\sqrt{(\ell - n)(\ell + n + 1)} \cdot t^{\ell}_{m, n+1}, \nonumber \\
	\partial_-(t^\ell_{mn})  = & -\sqrt{(\ell + n)(\ell - n + 1)} \cdot t^{\ell}_{m, n-1}, \nonumber \\
	\partial_0(t^\ell_{mn})  = & n \cdot t^\ell_{mn},
\end{align}
for all $\ell \in \frac{1}{2}\mathbb{N}_0$ and $m, n \in J_\ell$. 

As we will see below, the Lie group $SU(2)$ provides a rich setting for constructing examples of diagonal systems. Let us first discuss a class of systems of order $1$, inspired by Proposition \ref{ghgs-thm}.

Let $r,s \in \mathbb{N}$. Consider the compact Lie group $\mathcal{G}=\TrSs$ and the system $P = (P_1, P_2, \ldots, P_N)$, with
\begin{equation}\label{Lr-operators}
	P_\nu = \sum_{j=1}^r c_{\nu j} \partial_{x_j} + \sum_{k=1}^s id_{\nu k} \partial_{0,k} + q_\nu, \quad 1 \leq \nu \leq N,
\end{equation}
where the derivatives $\partial_{x_j} = \partial / \partial x_j$ are defined on distinct copies of the torus $\mathbb{T}^1$, the neutral operators $\partial_{0,k}$ are defined on distinct copies of the three-dimensional sphere $\mathbb{S}^3$. Moreover, the coefficients $c_{\nu j}, d_{\nu k}, q_\nu \in \mathbb{C}$, for $1 \leq j \leq r$, $1 \leq k \leq s$, and $1 \leq \nu \leq N$.

Let us translate the condition \eqref{DC2} in this case. Recall that we can identify $\widehat{\mathcal{G}}\simeq \mathbb{Z}^r \times (\frac{1}{2}\mathbb{N}_0)^s$. Fix $\Xi=(\xi,\ell)=(\xi_1,...,\xi_r,\ell_1,...,\ell_s) \in \mathbb{Z}^r \times (\frac{1}{2}\mathbb{N}_0)^s$. Since all elements in $\widehat{\mathbb{T}}$ have dimension $1$, we have $J_\Xi=\{(1,...,1,\alpha); \alpha \in J_\ell\}$, where $J_\ell = J_{\ell_1} \times... \times J_{\ell_s}$. 

If $e_{\xi}(x)=e^{ix \cdot \xi}$, then $\partial_{x_j} e_{\xi}\left. \right|_{x=0} = i\xi_j$ for each $j=1,...,r$. So $\mu_{(1,...,1)}(\xi)=\xi$.

Since $\sum_{k=1}^s id_{\nu k}\partial_{0,k} = \sum_{k=1}^s -d_{\nu k}D_{3,k}$ and by the relations (\ref{tabela}) we have $D_{3,k}(t^\ell_{\alpha_k \alpha_k})= -i \alpha_k t^\ell_{\alpha_k \alpha_k}$, for each $k=1,...,s$, then $\mu_\alpha(\ell) =- \alpha$. So
$$\mu_{(1,...,1,\alpha)}(\xi,\ell) = (\xi,-\alpha).$$

Denoting $c_\nu = (c_{\nu 1}, \ldots, c_{\nu r}) \in \mathbb{C}^r$ and $d_\nu = (d_{\nu 1}, \ldots, d_{\nu s}) \in \mathbb{C}^s$, we have
$$
\big|\jp{ (c_\nu, -d_\nu), \mu_{(1,...,1,\alpha)}(\xi,\ell)}-iq_\nu \big| = \big| \jp{\xi,c_\nu}+\jp{\alpha,d_\nu} - iq_\nu \big|,
$$ 
for all $(\xi,\ell) \in \mathbb{Z}^r \times (\frac{1}{2}\mathbb{N}_0)^s$, $\alpha \in J_\ell$ and $\nu=1,...,N$. 
In this way, the condition \eqref{DC2} takes the following form: 

there are constants $C>0$ and $\beta \in \mathbb{R}$ with the following property:
for all $(\xi,\ell) \in  \mathbb{Z}^r \times (\frac{1}{2}\mathbb{N}_0)^s$ such that $\max_\nu \big| \jp{\xi,c_\nu}+\jp{\alpha,d_\nu} - iq_\nu \big|>0$, there exists $\nu \in \{1,...,N\}$ such that
\begin{align}\label{toroxesfera}
	\big| \jp{\xi,c_\nu}+\jp{\alpha,d_\nu}  - iq_\nu \big| \geq C \cdot \left(\jp{\xi}+\jp{\ell} \right)^{-\beta}.
\end{align} 

The set $\mathcal{Z}$, defined in (\ref{zset}), can be rewritten as 
\begin{align}\label{zset_toroxesfera}
	\mathcal{Z}_{\bullet \bullet} = \left\{(\xi,\ell) %\in \mathbb{Z}^r \times (\tfrac{1}{2}\mathbb{N}_0)^s 
	\mid \exists \alpha \in J_\ell \textrm{ such that } \jp{\xi,c_\nu}+\jp{\alpha,d_\nu} - iq_\nu=0,  \nu=1,...,N \right\}.
\end{align}
This particular case has an important remark. If $\ell \in \frac{1}{2}\mathbb{N}_0$ and $k \in \mathbb{N}$, then $J_{\ell} \subset J_{\ell+k}$. So, in terms of vectors of indexes, for $\ell \in (\frac{1}{2}\mathbb{N})_0^s$, we have $J_{\ell} \subset J_{\ell+k}$, where $\ell+k\doteq \ell+(k,...,k)$. In particular, if $(\xi_0,\ell_0) \in \mathcal{Z}$, then $(\xi,\ell+k) \in \mathcal{Z}$ for all $k \in \mathbb{N}$. We conclude that if $\mathcal{Z} \neq \emptyset$, then $\mathcal{Z}$ is infinite. 

We summarize the discussion above in the following:

\begin{corollary}\label{cor_toroxesfera} For $N,r,s \in \mathbb{N}$ consider the operators
	$$
	P_\nu = \sum_{j=1}^r c_{\nu j} \partial_{x_j} + \sum_{k=1}^s id_{\nu k} \partial_{0,k} + q_\nu, \quad 1 \leq \nu \leq N,
	$$
	on the Lie group $\mathbb{T}^r \times \mathbb{S}^{3s}$ and the corresponding system $P=(P_1,...,P_N)$. Then:
	\begin{enumerate}
		\item $P$ is (GH) if, and only if, $P$ satisfies the condition (\ref{toroxesfera}) and the set $\mathcal{Z}_{\bullet \bullet}$ is empty;
		\item $P$ is (GS) if, and only if, $P$ satisfies the condition (\ref{toroxesfera}).
	\end{enumerate}
\end{corollary}

We now focus on examples to illustrate the conditions discussed in the previous results. The first two examples will apply the preceding corollary, while the third will showcase diagonal systems of arbitrary order.

First, let us recall the definition of Liouville numbers: we say that $\lambda \in \mathbb{R}$ is a non-Liouville number if there are positive constants $K$ and $C$ such that 
$$|a+\lambda \cdot b| \geq C |b|^{-K}$$
for all $a,b \in \mathbb{Z}$ with $b \neq 0$.  

\begin{example}
	For each $\nu \in \{1, 2, \ldots, N\}$, consider the operator  
	$$
	P_\nu = \alpha_\nu \partial_{x_\nu} +i \partial_0, \quad \text{where } \alpha_\nu \in \mathbb{R} \setminus \mathbb{Q},
	$$
	defined on the Lie group $\mathbb{T}^N \times \mathbb{S}^3$.  
	
	In the notation of Corollary \ref{cor_toroxesfera}, we have $ c_\nu = \alpha_\nu e_\nu $, $ d_\nu = 1 $, and $ q_\nu = 0 $ for all $\nu = 1, \ldots, N$. Consequently,  
	$$
	\mathcal{Z}_{\bullet \bullet} = \left\{ (\xi, \ell) \in \mathbb{Z}^N \times \tfrac{1}{2}\mathbb{N}_0 \; ; \; \exists m \in J_\ell \text{ such that } m + \xi_\nu \alpha_\nu = 0, \nu =1, \ldots, N \right\}.
	$$  
	
	Since $\alpha_\nu \in \mathbb{R} \setminus \mathbb{Q}$, it is clear that $ m + \xi_\nu \alpha_\nu = 0 $ for all $\nu = 1, \ldots, N$ if, and only if, $\xi = 0$ and $m = 0$. Therefore, $(0, \ell) \in \mathcal{Z}_{\bullet \bullet}$ for all $\ell \in \mathbb{N}$. By Corollary \ref{cor_toroxesfera}, this shows that the system $P = (P_1, \ldots, P_N)$ is not (GH), regardless of the specific values of the irrational numbers $\alpha_\nu$.  
	
	On the other hand, the situation changes when studying the global solvability of $P$. Assume now that each $\alpha_\nu$ is not a Liouville number for $\nu = 1, \ldots, N$. In this case, there exist positive constants $C_\nu$ and $K_\nu$ such that  
	$$
	|a + \alpha_\nu b| \geq C_\nu |b|^{-K_\nu}
	$$  
	for all $a, b \in \mathbb{Z}$ with $b \neq 0$.  
	
	Now fix $(\xi, \ell) \in \mathbb{Z}^N \times \frac{1}{2}\mathbb{N}_0$, $m \in J_\ell$, and $\nu \in \{1, \ldots, N\}$ such that $m + \xi_\nu \alpha_\nu \neq 0$. In this scenario, we must analyze the following cases:  
	\begin{enumerate}
		\item $\xi_\nu=0$ and $m \neq 0$. In this case
		$$|m+\xi_\nu \alpha_\nu| = |m| \geq \frac{1}{2};$$
		\item $\xi_\nu \neq 0$ and $m=0$. In this case
		$$|m+\xi_\nu \alpha_\nu| = |\xi_\nu | \cdot |\alpha_\nu| \geq |\alpha_\nu|;$$
		\item $\xi_\nu \neq 0$ and $m\neq 0$. In this case
		$$
		|m+\xi_\nu \alpha_\nu| = \frac{1}{2} | (2m)+(2\xi_\nu)\alpha_\nu| \geq 
		\frac{C_\nu}{2^{K_\nu+1}} |\xi_\nu|^{-K_\nu} \geq  \frac{\widetilde{C}_\nu}{2^{K_\nu+1}}(\jp{\xi}+\jp{\ell})^{-K_\nu}.
		$$
	\end{enumerate}
	
	Defining $C\doteq \displaystyle\min_{1 \leq \nu \leq N} \left\{\frac{1}{2},|\alpha_\nu|, \widetilde{C}_\nu \cdot 2^{-K_\nu-1}\right\}>0$ and $\beta = \displaystyle\max_{1 \leq \nu \leq N} K_\nu>0$, then in any of the three cases above we obtain
	$$|m+\xi_\nu \alpha_\nu| \geq C \left(\jp{\xi}+\jp{\ell} \right)^{-\beta}.$$
	proving that condition \eqref{DC2} holds and by Corollary \ref{cor_toroxesfera} the system $P=(P_1,...,P_N)$ is (GS). 
\end{example} 	

\begin{example}
	In the notation of Corollary \ref{cor_toroxesfera}, let us consider the case $r = 2$, $s = 1$, and the system of operators  
	$$
	\begin{cases}
		P_1 & = \partial_{x_1} + i\alpha_1 \partial_0 + q_1,\\
		P_2 & = \partial_{x_2} + i\alpha_2 \partial_0 + q_2,
	\end{cases}
	$$
	where $\alpha_1, \alpha_2, q_1, q_2$ are complex constants. For $\nu = 1, 2$, $(\xi, \ell) \in \mathbb{Z}^2 \times \frac{1}{2}\mathbb{N}_0$, and $m \in J_\ell$, we have  
	$$
	|\jp{\xi, c_\nu}+\jp{\alpha, d_\nu} - i q_\nu| = |\big(m \operatorname{Re}(\alpha_\nu) + \operatorname{Im}(q_\nu) - \xi_\nu\big) 
	+ i\big(m \operatorname{Im}(\alpha_\nu) - \operatorname{Re}(q_\nu)\big)|.
	$$
	
	Now assume that $\operatorname{Im}(\alpha_\nu) \neq 0$ and $\operatorname{Re}(q_\nu) / \operatorname{Im}(\alpha_\nu) \notin \frac{1}{2} \mathbb{Z}$ for $\nu = 1, 2$. Under these assumptions, we conclude that  
	$$
	\inf_{m \in \frac{1}{2}\mathbb{Z}} |m - \operatorname{Re}(q_\nu) / \operatorname{Im}(\alpha_\nu)| \doteq K_\nu > 0.
	$$  
	
	Letting $K = \min\{K_1, K_2\}$ and $C = \min\{|\operatorname{Im}(\alpha_1)|, |\operatorname{Im}(\alpha_2)|\}$, it follows that  
	$$
	|\jp{\xi, c_\nu}+\jp{\alpha, d_\nu}  - i q_\nu| \geq |\operatorname{Im}(\alpha_\nu)| \cdot |m - \operatorname{Re}(q_\nu) / \operatorname{Im}(\alpha_\nu)| 
	\geq |\operatorname{Im}(\alpha_\nu)| \cdot K_\nu 
	\geq C K > 0,
	$$  
	for all $\nu = 1, 2$ and $(\xi_1, \xi_2, m) \in \mathbb{Z} \times \mathbb{Z} \times \frac{1}{2} \mathbb{Z}$.  
	
	Therefore, in this case, $\mathcal{Z}_{\bullet \bullet} = \emptyset$, and condition \eqref{DC2} is satisfied. By Corollary \ref{cor_toroxesfera}, the system $P$ is (GH).  
\end{example}

\begin{example} Here we present a class of diagonal system of any order. Let $p,r \in \mathbb{N}$ and $a_1,...,a_r,b_1,...,b_r$ be non-zero real numbers. For each $j=1,...,r$ consider the operator $P_j: \mathcal{D}'(\mathbb{S}^3)\rightarrow \mathcal{D}'(\mathbb{S}^3)$ defined by 
	$$P_j = a_j\partial_0^p + ib_j \partial_+ \partial_-.$$
	and the corresponding system $P=(P_1,...,P_r)$. By the relations (\ref{tabela}) we get for each $j=1,...,r$ and $\ell \in \frac{1}{2}\mathbb{N}_0$ that the matrix $\sigma_{P_j}(\ell)$ is diagonal  and for each $n \in J_\ell$ we have
	$$\sigma_n^j(\ell) = a_jn^p +ib_j (\ell+n)(\ell-n+1).$$
	Since $a_j,b_j$ are non-zero real numbers, for $\ell \in \frac{1}{2}\mathbb{N}$ and $n\neq 0$ we have
	$$|\sigma_n^j(\ell)| \geq |a_j| |n|^p \geq |a_j| |1/2|^p \geq a/2^p >0$$
	where $a=\min_j |a_j|>0$. 
	For $n=0$ we have 
	$$|\sigma_n^j(\ell)| = |b_j|\cdot \ell \cdot (\ell+1) > |b_j|\geq b$$
	where $b=\min_j |b_j|>0$. So, if $C=\min \{a/2^p,b\}>0$, then for all $\ell \in \frac{1}{2}\mathbb{N}$ and $n \in J_\ell$ we have
	$$|\sigma_n^j(\ell)| \geq C > C \jp{\ell}^{-1}.$$
	In particular, $\mathcal{Z}$ is finite and the system $P$ satisfies the condition \eqref{DC1}, so by the Theorem \ref{ghthm} the system $P$ is (GH).
	
\end{example}

%==========================================================	
%==========================================================	
\section{Global Properties for Triangular Systems}
%==========================================================	
%==========================================================	

In this section, we consider triangular systems of left-invariant operators defined on compact Lie groups under different additional hypothesis about the group or about the control of the symbol. With these requirements, we establish sufficient conditions to ensure the global hypoellipticity of the system.

\begin{definition}
	A system of left-invariant continuous operators $P = (P_1, \dots, P_r)$ is said to be triangular if, for each $\Xi \in \widehat{G}$, there exists a matrix representative $\xi \in \Xi$ such that the matrix $\sigma_{P_j}(\xi)$ is triangular for all $j = 1, \dots, r$. In this case, for each $\Xi = [\xi] \in \widehat{G}$ and $1 \leq j \leq r$, we write  
	$$ 
	\sigma_{P_j}(\xi) = \left( \sigma_{k\ell}^j(\xi)\right)_{1 \leq k, \ell \leq d_\xi}.
	$$ 
\end{definition}

The diagonal entries are denoted as in previous sections, that is, $ \sigma_k^j(\xi) \doteq \sigma_{kk}^j(\xi) $.  

\begin{example}
	As established in \cite[Section 1.1]{RRose-Simult-triang_book} and \cite[Section 6.5]{Hoffman-Kunza-book}, if the family of left-invariant operators $P_1, \dots, P_r$ commutes, then $P = (P_1, \dots, P_r)$ is a triangular system. In particular, every diagonal system is also a triangular system.
\end{example}

%	\begin{definition} 
	%		We say that a triangular system $P = (P_1, \dots, P_r)$ satisfies the $(DC_T)$ condition if there exist positive constants $C$ and $N$ such that  
	%		\begin{equation}\label{DC_T} 
		%			\big|\sigma_k^j(\xi) \big| \geq C \jp{\xi}^{-N} \tag{DC$_T$}
		%		\end{equation} 
	%		for all $[\xi] \in \widehat{G}$, $1 \leq k \leq d_\xi$, and $j = 1, \dots, r$ such that $\sigma_k^j(\xi) \neq 0$.  
	%	\end{definition}
%
%Note that condition \eqref{DC_T} reduces to condition \eqref{DC1}, as stated in Definition \ref{DC-diagonal-def}, when restricted to diagonal systems. In other words, as expected, the asymptotic behavior of solutions to triangular systems is determined solely by the eigenvalues of the matrix.

%\begin{definition} \label{C-and-N_sigma}
%	We say that triangular system $P = (P_1, \dots, P_r)$ satisfies the condition $(DC_T)$, if there are constants $C>0$ and $M \in \mathbb{R}$ such that  
%	\begin{equation*}
%		\big|\sigma_k^j(\xi) \big| \geq C \jp{\xi}^{-M} 
%	\end{equation*} 
%	for all $[\xi] \in \widehat{G}$, $1 \leq k \leq d_\xi$, and $j = 1, \dots, r$ such that $\sigma_k^j(\xi) \neq 0$.  
%\end{definition}
	
	If $m_j$ is the order of each operator $P_j$, by setting $\widetilde{M} = \max \{m_1, \dots, m_r\}$, there exists a constant $\widetilde{C} \geq 0$ such that  
	$$
	\big| \sigma_{k\ell}^j(\xi) \big| \leq \widetilde{C} \jp{\xi}^{\widetilde{M}}
	$$
	for all $1 \leq k, \ell \leq d_\xi$, $j = 1, \dots, r$, and $\Xi = [\xi] \in \widehat{G}$.  
	
	If the system $P$ satisfies the condition (\ref{DC1}) then there exists $M \in \mathbb{R}$ and $C>0$ such that
	$$\max_{j\in \bigwedge_{k,\xi}} |\sigma_k^j(\xi)| \geq C \cdot \jp{\xi}^{-M}$$
	whenever $\bigwedge_{k,\xi} \neq \emptyset$.  
	
	Thus, if we set $C_\sigma \doteq \max \{1/C, \widetilde{C}\}$ and $N_\sigma = \max \{M, \widetilde{M}\}$, then  
	\begin{equation}\label{dct1}
		\big|\sigma_{k\ell}^j(\xi)\big| \leq C_\sigma \jp{\xi}^{N_\sigma}
		\quad \text{and} \quad			
		\frac{1}{\max_{j\in \bigwedge_{k,\xi}} |\sigma_k^j(\xi)|} \leq C_\sigma \jp{\xi}^{ N_\sigma},
	\end{equation}
	for all $1 \leq k, \ell \leq d_\xi$, $j = 1, \dots, r$, and $\Xi = [\xi] \in \widehat{G}$ for the first inequality, and the second one holds whenever $\bigwedge_{k,\xi} \neq \emptyset$. 
	
	These estimates will be used in the proof of the main result of this section.

Now, we introduce an additional hypothesis related to the boundedness of the dimensions of representations.

\begin{definition}\label{B-condition} 
	A compact Lie group $G$ satisfies the bounded dimension condition \eqref{BD} if there exists a constant $d_G > 0$ such that 
	\begin{equation}\label{BD}
		d_\xi \leq d_G, \quad \text{for all } \Xi = [\xi] \in \widehat{G}. \tag{BD}
	\end{equation} 
\end{definition} 

This class of groups was fully characterized by C. Moore in \cite{moore1972groups}. Specifically, a compact Lie group $G$ satisfies the condition \eqref{BD} if, and only if, $G$ contains an open abelian subgroup of finite index. For instance, this class includes any semi-direct product $A \rtimes H$, where $A$ is an abelian compact Lie group and $H$ is a finite group.

Finally, to state the main result of this section, we recall that the set \eqref{zset}, introduced in Theorem \ref{Zinfinite-notGH}, is defined in this case as  
\begin{equation*}\label{zset2} 
	\mathcal{Z} \doteq \{ [\xi] \in \widehat{G} \mid \exists k_0 \in \{1, \dots, d_\xi\}, \ \text{such that } \sigma_{k_0}^j(\xi) = 0 \ \text{for all } j = 1, \dots, r \},
\end{equation*}
which coincides with the set $\mathcal{Z}$ in the case of diagonal systems.

\begin{theorem}\label{GH-triang-system}
	Let $G$ be a compact Lie group that satisfies the condition \eqref{BD}. If $P = (P_1, \dots, P_r)$ is a triangular system that satisfies the condition \eqref{DC1} and that the set $\mathcal{Z}$ is finite,  then $P$ is (GH).
\end{theorem}

\begin{proof}
	Without loss of generality, we may assume that the matrices $\sigma_{P_j}(\xi)$ are lower-triangular for all $[\xi] \in \widehat{G}$ and $j = 1, \dots, r$. 
	
	Let $f = (f_1, \dots, f_r) \in C^\infty(G)^r$ and $u \in \mathcal{D}'(G)$ such that $P u = f$. To prove that $u$ is actually a smooth function, it is enough to verify that its Fourier coefficients satisfy the smoothness condition \eqref{smooth1}. Specifically, for each $N > 0$, there is a constant $C_{N} > 0$, depending only on $N$, such that
	\begin{equation}\label{u-desired-estimate}
		\big| \widehat{u}(\xi)_{k\ell} \big| \leq C_{N} \langle \xi \rangle^{-N},		
	\end{equation}
	for all $[\xi] \in \widehat{G}$ and $1 \leq k, \ell \leq d_\xi$.
	
	First, since $f \in C^\infty(G)^r$, for any $N >0$, there exists a constant $C_{N,f} > 0$, depending only on $N$, such that 
	\begin{equation}\label{cond_f}
		\left| \widehat{f_j}(\xi)_{k\ell} \right| \leq C_{N,f} \langle \xi \rangle^{-N}
	\end{equation}
	for all $[\xi] \in \widehat{G}$, $1 \leq k, \ell \leq d_\xi$, and $j = 1, \dots, r$. 
	
	We now proceed to estimate $\big| \widehat{u}(\xi)_{k\ell} \big|$ recursively. Starting with the case $k = 1$, for each $[\xi] \notin \mathcal{Z}$ we have $\bigwedge_{1,\xi} \neq \emptyset$, so there exists $j(1) \in \{1, \dots, r\}$ such that 
	$$
	\sigma_1^{j(1)}(\xi) \neq 0 \textrm{ and } \frac{1}{|\sigma_1^{j(1)}(\xi)|} \leq C_\sigma \jp{\xi}^{N_\sigma} 
	$$
	
	Using the triangular structure of $P$ and the fact that $Pu = f$, the Fourier coefficient $\widehat{u}(\xi)_{1\ell}$ satisfies:
	$$
	\sigma_1^{j(1)}(\xi) \cdot \widehat{u}(\xi)_{1\ell} = \widehat{f}_{j(1)}(\xi)_{1\ell}.
	$$
	
	From (\ref{dct1}) and (\ref{cond_f}), we deduce
	\begin{align}
		\big| \widehat{u}(\xi)_{1\ell} \big| & \leq  \frac{1}{\big| \sigma_1^{j(1)}(\xi) \big|} \cdot \big| \widehat{f}_{j(1)}(\xi)_{1\ell} \big|  \label{u1L} \\
		&  \leq  C_\sigma \langle \xi \rangle^{N_\sigma} \cdot C_{N,f} \langle \xi \rangle^{-N} \nonumber \\
		& =  C_{N}^1 \langle \xi \rangle^{-N+N_\sigma}, \nonumber
	\end{align}
	where $C_N^1 \doteq C_\sigma C_{N,f}$ depends only on $N$.
	
	Next, consider the case $k = 2$. Since $[\xi] \notin \mathcal{Z}$, we have $\bigwedge_{2,\xi} \neq \emptyset$ so there exists $j(2) \in \{1, \dots, r\}$ such that 
	$$
	\sigma_2^{j(2)}(\xi) \neq 0 \textrm{ and } \frac{1}{|\sigma_2^{j(2)}(\xi)|}\leq C_\sigma \jp{\xi}^{N_\sigma}
	$$  
	From the triangular structure of $P u = f$, we have the equation
	$$
	\sigma_{2, 1}^{j(2)}(\xi) \cdot \widehat{u}(\xi)_{1\ell} + \sigma_2^{j(2)}(\xi) \cdot \widehat{u}(\xi)_{2\ell} = \widehat{f}_{j(2)}(\xi)_{2\ell}.
	$$	
	Applying the estimates from (\ref{dct1}), (\ref{cond_f}), and \eqref{u1L}, we obtain
	\begin{align}
		\big| \widehat{u}(\xi)_{2\ell} \big| & \leq \frac{1}{\big| \sigma_2^{j(2)}(\xi) \big| } 
		\left( \big| \widehat{f}_{j(2)}(\xi)_{2\ell} \big| + \big| \sigma_{2, 1}^{j(2)}(\xi) \big| \cdot \big|\widehat{u}(\xi)_{1\ell} \big| \right) \\
		& \leq C_\sigma \langle \xi \rangle^{N_\sigma} \left(  C_{N,f} \langle \xi \rangle^{-N} 
		+ C_\sigma \langle \xi \rangle^{N_\sigma} \cdot C_{N}^1 \langle \xi \rangle^{-N+N_\sigma} \right) \nonumber \\
		& \leq C_\sigma C_{N,f}  \langle \xi \rangle^{-N+N_\sigma} 
		+ (C_\sigma)^2 C_{N}^1 \langle \xi \rangle^{-N+3N_\sigma} \nonumber \\
		& \leq \left(C_\sigma C_{N,f} + (C_\sigma)^2 C_{N}^1\right) \langle \xi \rangle^{-N+3N_\sigma} \nonumber \\ 
		&= C_N^2 \langle \xi \rangle^{-N+3N_\sigma}, \nonumber 		 
	\end{align}
	where $C_N^2 \doteq \left(C_\sigma C_{N,f}  + (C_\sigma)^2 C_{N}^1\right)$ depends only on $N$.
	
	Using the \eqref{BD} condition, this process can be repeated up to $d_G - 2$ additional steps, resulting in the general estimate:
	$$
	\big| \widehat{u}(\xi)_{k\ell} \big| \leq C_N^k \langle \xi \rangle^{-N+(2k-1) N_\sigma},
	$$
	for all $1 \leq k, \ell \leq d_\xi$ and $[\xi] \notin \mathcal{Z}$.
	
	Therefore, for all $1 \leq k, \ell \leq d_\xi$ and $[\xi] \notin \mathcal{Z}$, we have
	$$
	\big| \widehat{u}(\xi)_{k\ell} \big| \leq C_{N,d_G} \langle \xi \rangle^{-N+(2 d_G-1) N_\sigma},
	$$
	where $C_{N,d_G} \doteq \max \big\{ C_N^k \mid 1 \leq k \leq d_G \big\} > 0$ is a constant that depends only on $N$.
	
	To conclude the proof, note that since the set $\mathcal{Z}$ is finite, we can define 
	$$
	C_\mathcal{Z} \doteq \max \left\{ \left| \widehat{u}(\xi)_{k\ell} \right| \mid  [\xi] \in \mathcal{Z}, \ 1 \leq k, \ell \leq d_\xi \right\} \geq 0.
	$$
	
	By setting $C_N \doteq \max \{C_{N,d_G}, C_\mathcal{Z} \}$, it follows that 
	$$
	\big| \widehat{u}(\xi)_{k\ell} \big| \leq C_{N} \langle \xi \rangle^{-N+(2 d_G-1) N_\sigma},
	$$
	for all $[\xi] \in \widehat{G}$ and $1 \leq k, \ell \leq d_\xi$. 
	
	Since $C_N$ depends only on $N$, and $(2d_G - 1) \geq 0$ is fixed, this estimate ensures that $u \in C^\infty(G)$, thereby completing the proof.	
\end{proof}

Next, we turn our attention to the range of $P$. The idea is to prove that $P$ is almost-globally hypoelliptic, as defined in Remark \ref{AGH->GS}.

\begin{theorem}\label{almostGH-triang-system}
	Let $G$ be a compact Lie group that satisfies the condition \eqref{BD}. If a triangular system  $P = (P_1, \dots, P_r)$ satisfies the condition \eqref{DC1}, then $P$ is almost-globally hypoelliptic.
\end{theorem}

Note that the assumption of the finiteness of the set $\mathcal{Z}$, required in the previous theorem, is essential to ensure global hypoellipticity. Consequently, the existence of distributional solutions that are not smooth is related to the Fourier frequencies arising from representations associated with this set. In this context, we consider the following set:
$$
\ker \mathbb{P}(\xi) = \bigcap_{j=1}^r \mathrm{Ker} (\sigma_{P_j} (\xi)) 
$$
and its orthogonal complement $ \ker \mathbb{P}(\xi)^\perp $.

\begin{proof}
	Let $ u \in \mathcal{D}'(G) $ such that $ P u = f $. We consider the following decomposition of its Fourier coefficients:
	$$
	\widehat{u}(\xi)_{k\ell} = \widehat{v}(\xi)_{k\ell} + \widehat{w}(\xi)_{k\ell},
	$$
	with $ \widehat{v}(\xi)_{k\ell} \in \ker \mathbb{P}(\xi)^\perp $ and $ \widehat{w}(\xi)_{k\ell} \in \ker \mathbb{P}(\xi) $, for all $ [\xi] \in \widehat{G} $ and $ 1 \leq k, \ell \leq d_\xi $.
	
	Note that the distribution $ v \in \mathcal{D}'(G) $, whose Fourier coefficients $ \widehat{v}(\xi) $ are given above, satisfies $ P_j u = P_j v $ for all $ j = 1, \dots, r $.  Thus, it suffices to repeat the arguments of the proof of the previous theorem to conclude that the distribution $ v \in C^\infty(G) $.	
\end{proof}

As a consequence of this theorem and Remark \ref{AGH->GS}, we obtain the following result:

\begin{corollary}
	Let $G$ be a compact Lie group satisfying the condition \eqref{BD}. If a triangular system $P = (P_1, \dots, P_r)$ satisfies the condition \eqref{DC1}, then $P$ has a closed range.
\end{corollary}

Next, we will present another Theorem about global hypoellipticity for the triangular system $P$ without assuming the condition \eqref{BD}.

Suppose $P$ is a triangular system on a compact Lie group $G$ such that $\mathcal{Z}$ is finite and condition (\ref{DC1}) holds with constants $C>0$ and $M>0$. In this case, for all $[\xi] \in \widehat{G}$ such that $\bigwedge_{k,\xi} \neq \emptyset$, we have 
$$ \frac{1}{|\sigma_k^j(\xi)|} \leq \frac{1}{C\jp{\xi}^M}.$$
for some $j \in \{1,...,r\}$. 
Since $M>0$ and $\jp{\xi}>1$, there exist $R>0$ such that $\jp{\xi}>R$ implies $1/(C\jp{\xi}^{M})<1$. So, for $\jp\{\xi\}>R$ and $k$ such that $\bigwedge_{k,\xi} \neq \emptyset$ we can assume there exists $j \in \{1,...,r\}$ such that $1/|\sigma_k^j(\xi)| <1$.  

Now, suppose that the nilpotent part of $\sigma_{P}$ is uniformely bounded, that is, there exists $K>0$ such that $|\sigma_{k\ell}^j(\xi)| \leq K$ for all $[\xi] \in \widehat{G}$, $1 \leq k,\ell \leq d_\xi$ with $k>l$. It is clear that $P$ is $(GH)$ if, and only if, $\widetilde{P}:=(1/K) \cdot P$ is $(GH)$. So, we may assume that $K=1$. 

Suppose that $f=(f_1,...,f_r) \in C^\infty(G)^r$ and $u \in \mathcal{D}'(G)$ are such that $Pu=f$ and fix $N>0$. Since $f \in C^\infty(G)^r$, there exists $C(N)>0$ such that $|\widehat{f}(\xi)_{k\ell}| \leq C(N) \jp{\xi}^{-N}$ for all $[\xi] \in \widehat{G}$, $1\leq k,\ell \leq d_\xi$ and $j=1,...,r$. Normalizing again, there exists $\widetilde{R}>0$ such that $|\widehat{f}(\xi)_{k\ell}| \leq \jp{\xi}^{-N+1}$ for all $[\xi]$ with $\jp{\xi}> \widetilde{R}$, $1 \leq k,\ell \leq d_\xi$ and $j=1,...,r$. Define $R_N:=\max \{R,\widetilde{R}\}$. 

Like in the previous proofs, it is enought to estimate $|\widehat{u}(\xi)_{k\ell}|$ for $[\xi] \in \widehat{G}\setminus \mathcal{Z}$. So let $[\xi] \in \widehat{G}\setminus \mathcal{Z}$ such that $\jp{\xi}>R_N$ and fix $1 \leq \ell \leq d_\xi$. Since $[\xi] \notin \mathcal{Z}$, we have $\bigwedge_{1,\xi} \neq \emptyset$, so there exists $j$ such that $\sigma_1^j(\xi) \neq 0$ and $1/|\sigma_1^j(\xi)|<1$. Since $\sigma_{P_{j_1}}$ is lower triangular, we have
$$\sigma_1^j(\xi) \cdot \widehat{u}_{1\ell}(\xi) = \widehat{f}_{j}(\xi)_{1\ell}.$$
So 
$$|\widehat{u}_{1\ell}(\xi)| \leq \frac{1}{|\sigma_1^j(\xi)|} \cdot |\widehat{f}_{j}(\xi)_{1\ell}| <\jp{\xi}^{-N+1}.$$
Now for $k=2$, again since $[\xi] \in \widehat{G}\setminus \mathcal{Z}$, we have $\bigwedge_{2,\xi} \neq \emptyset$ and there exists $j$ with $\sigma_2^{j}(\xi)\neq 0$ and $1/|\sigma_2^j(\xi)|<1$. So
$$\sigma_{21}^{j}\widehat{u}(\xi)_{1\ell}+\sigma_2^{j}(\xi) \cdot \widehat{u}(\xi)_{2\ell} = \widehat{f}_{j}(\xi)_{2\ell}$$
which implies
$$|\widehat{u}(\xi)_{2\ell}| \leq \frac{1}{|\sigma_2^{j}(\xi)|}\cdot \left[|\widehat{f}_{j}(\xi)_{2\ell}|+|\sigma_{21}^{j}|\cdot |\widehat{u}(\xi)_{1\ell}| \right] \leq 1 \cdot \left[  \jp{\xi}^{-N+1}+1 \cdot\jp{\xi}^{-N+1}  \right]=2 \cdot \jp{\xi}^{-N+1}.$$
In general, it is clear that
$$|\widehat{u}(\xi)_{k\ell}| \leq k \cdot \jp{\xi}^{-N+1} \leq d_\xi \cdot \jp{\xi}^{-N+1}.$$
By Proposition 10.3.19 from \cite{RT2010_book}, there exists $C>0$ such that $d_\xi \leq C \jp{\xi}^{\dim(G)/2}$ for all $[\xi] \in \widehat{G}$, so
$$|\widehat{u}(\xi)_{k\ell}| \leq C \jp{\xi}^{-N+1+\dim(G)/2}$$
for all $[\xi] \notin \mathcal{Z}$ with $\jp{\xi}>R_N$, $1 \leq k,l \leq d_\xi$. This implies that $u \in C^\infty(G)$ and so $P$ is $(GH)$.

The previous argument proved the following:

\begin{theorem} Let $G$ be a compact Lie group and $P: \mathcal{D}'(G) \rightarrow \mathcal{D}'(G)^r$ be a continuous left-invariant triangular system on $G$. If
	\begin{itemize}
		\item[1)] the set $\mathcal{Z}$ is finite;
		\item[2)] the condition (\ref{DC1}) holds;
		\item[3)] the nilpotent part of $\sigma_P$ is uniformely bounded;	
	\end{itemize}
then $P$ is globally hypoelliptic.	
\end{theorem}

%============================================================
%============================================================
\section*{Declarations}
%============================================================
%============================================================

\subsection*{Funding}
The first author was supported in part by FAPESP (grant 2018/14316-3) and CNPq (grant 313581/2021-5). The second author was supported in part by CNPq (grants 316850/2021-7 and 423458/2021-3). The third author was supported in part by CNPq (grant 423458/2021-3) and a postdoctoral fellowship (grant 150776/2024-1).

\subsection*{Availability of data and material}
Not applicable.

\subsection*{Competing interests}
The authors declare that they have no competing interests.

\subsection*{Contributions of the authors}
All authors contributed equally to this work. All authors read and approved the final manuscript.

\bibliography{references}

\end{document}